\documentclass{amsart}

\usepackage{amsmath}
\usepackage{amssymb}
\usepackage{graphicx}
\usepackage{url}

\newtheorem{theorem}{Theorem}[section]
\newtheorem{lemma}[theorem]{Lemma}
\newtheorem{corollary}[theorem]{Corollary}

\theoremstyle{definition}

\theoremstyle{remark}
\newtheorem{remark}[theorem]{Remark}

\numberwithin{equation}{section}

\newfont{\hiera}{cmsy10 scaled 2488}
\newfont{\hierb}{cmsy10 scaled 1728}
\newfont{\hierc}{cmsy10 scaled 1200}

\newcommand{\Bigast}{
\mathop{\vphantom{\sum}\lower2.5pt\hbox{\hiera\char3}}}%

\newcommand{\Bigtimes}{
\mathop{\vphantom{\sum}\lower2.5pt\hbox{\hiera\char2}}}%

\begin{document}




\title[Minimal Polynomials for the Harborth Graph]{Minimal Polynomials for the Coordinates of the Harborth Graph}

\author{Eberhard H.-A. Gerbracht}
\curraddr{Bismarckstra\ss e 20, D-38518 Gifhorn, Germany}
\email{e.gerbracht@web.de}


\date{}


\begin{abstract}
The Harborth graph is the smallest known example of a 4-regular planar unit-distance graph. In this paper we give an analytical description of the coordinates of its vertices for a particular embedding in the Euclidean plane. More precisely, we show how to calculate the minimal polynomials of the coordinates of its vertices with the help of a computer algebra system, and list these. Furthermore some algebraic properties of these polynomials, and consequences to the structure of the Harborth graph are determined.
\end{abstract}



\maketitle

\section{Introduction}
The Harborth graph (see Figure 1) is the smallest known example of a 4-regular planar unit-distance graph. That is a planar graph, all of which edges are of unit length, with exactly four edges meeting in each vertex. This graph was named after its discoverer H.\ Harborth, who first presented it to the general public as a research problem in \cite{HarbStreichh}  and to a large international audience in a talk at the Eug\`ene Strens Memorial Conference on Recreational Mathematics and its History in 1986 \cite{HarbMatch,MathLand}. At both occasions he posed the question, if a smaller example of a 4-regular planar unit-distance graph could be found.

\begin{figure}
\centering
\includegraphics[width=.8\linewidth]{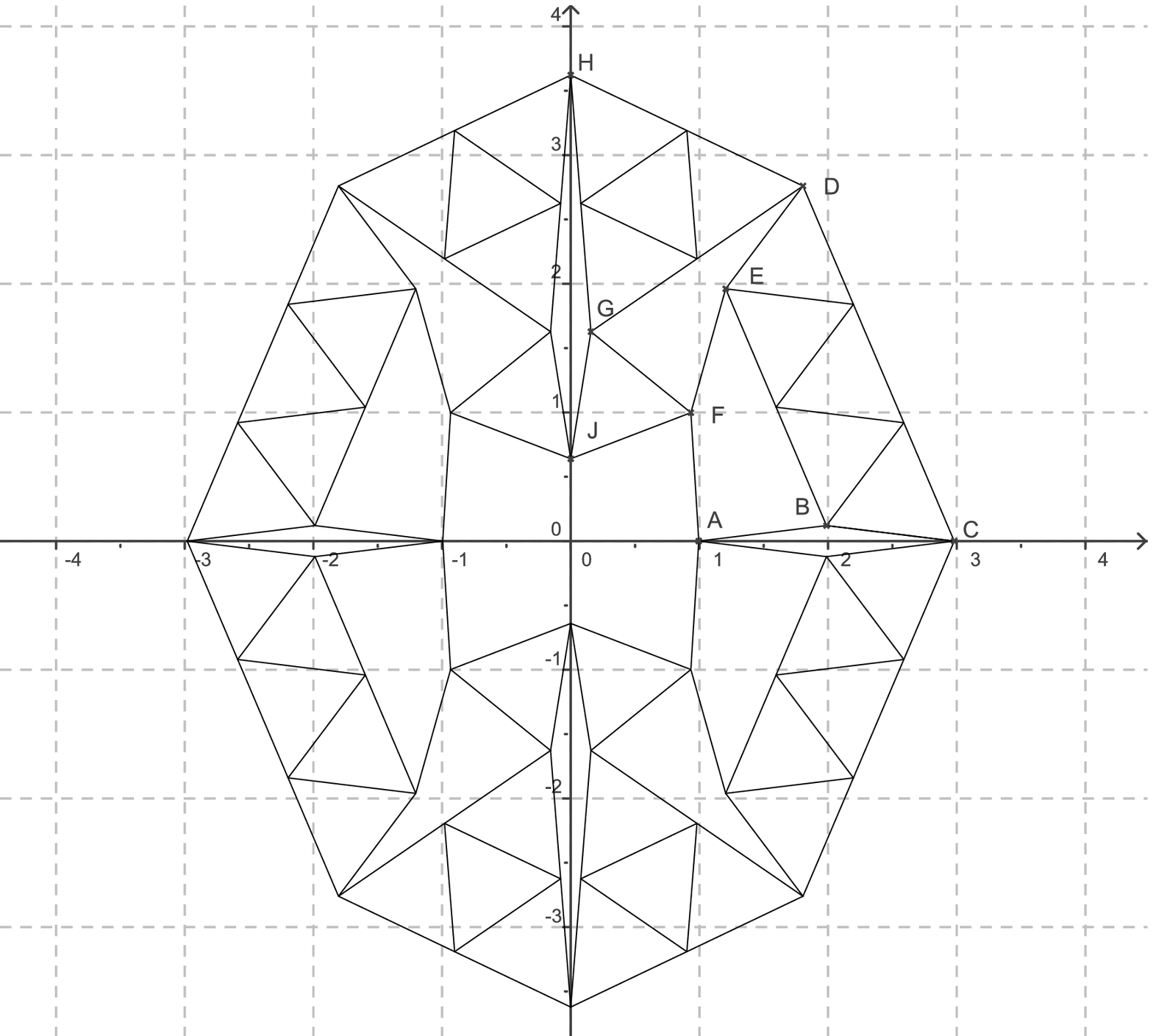}
\caption{The Harborth graph embedded in the Euclidean plane.}
\label{HarborthComplete}
\end{figure}

Curiously enough, up until now nearly all published pictures of the Harborth graph  -- even the original ones in \cite{HarbStreichh,HarbMatch}, as well as those in textbooks \cite{HarbKem, Pearls} -- seem to be slightly vague and inaccurate, with the vertices always being depicted by large dots. Furthermore there has not been given any analytical description of the Harborth graph, yet. That is, if we consider the graph as embedded in the Euclidean plane with a given coordinate system, the coordinates of the vertices have never been calculated exactly. This has gone to the point that the Harborth graph was even thought by some to be nonrigid\footnote{
If we consider the Harborth graph as a (mechanical) framework consisting of rigid bars interconnected by rotable joints, nonrigidity means that some vertices can be moved with respect to each other, that is the whole framework allows motions which are different from congruences  \cite{CombEEMech}.}, which, as the results of this paper imply,  cannot be the case.
Thus the wish for an ''exact'' description remained. This wish has recently been expressed again in the world wide web \cite{MathPuzzle}, which prompted this paper.

With the advent of dynamical geometry systems, several authors were finally able to produce more precise pictures \cite{HarbGeoExp,MathPuzzle}, which led to further evidence that there is one unique realisation of the Harborth graph in the Euclidean plane. With this paper we go one step further: using one particular way of construction, we set up a set of quadratic equations which completely describe the coordinates of crucial vertices of the Harborth graph. Using this initial set of equations, we will show that all the coordinates are algebraic numbers, and we will calculate their minimal polynomials with the help of a computer algebra system. Even though, as we will see later on, it is impossible to solve the corresponding algebraic equations exactly, which in our understanding means in terms of radicals, nevertheless we thus have achieved an exact analytic description of the Harborth graph, since, together with easy to calculate numerical approximations for the actual coordinates, these polynomials uniquely determine each coordinate.

The author is indebted to C.\ Adelmann and H.\ L\"owe from the Technical University Braunschweig, Germany, for a number of very valuable discussions on the subject. Furthermore he owes thanks to the Institute for Mathematical Physics of the TU Braunschweig, especially to the research group of R.F.\ Werner, for making its computing facilities available to him for this research.
  
\section{Using Dynamic Geometry Software and Numerical Analysis}

\subsection{Geometric construction}
Because of the obvious twofold symmetry of the Harborth graph, it is enough to analyse one of its quarters. Therefore, in a first step we construct one of these quarters (see Figure 2), using one of the existing (imprecise) first generation pictures 
as a blueprint.

 We start from an initial isosceles triangle $ABC$ of fixed but arbitrary height $T,$ with two sides being of unit length, and a neighbouring symmetric trapezoid $BCDE,$ which has the side $BC$ with the initial triangle in common. The parallel sides of the trapezoid are chosen to be of length $2$ and $3$ respectively. The remaining points are constructed from this initial configuration by using compass and ruler techniques. In the following, we list the necessary steps. Thereby $\hbox{Circ}(P,r)$ denotes the circle with center $P$ and radius $r,$ and $\cap$ the operation of letting two geometric figures intersect. Thus we get
\begin{eqnarray}
F &:=& \hbox{Circ}(A,1) \cap \hbox{Circ}(E,1) \label{Fconst}\\
G &:=& \hbox{Circ}(F,1) \cap \hbox{Circ}(D,2) \label{Gconst}\\
H &:=& \hbox{Circ}(D,2) \cap \hbox{Circ}(G,2) \label{Hconst}\\
J &:=& \hbox{Circ}(F,1) \cap \hbox{Circ}(G,1) \label{Jconst}
\end{eqnarray}

Although the intersection of two circles usually consists of two points, we use this notation as if there was no ambiguity. We are allowed to this, because we choose the resulting points of intersection according to our blueprint. Thus, e.g., $F$ is chosen in such a way that the quadrangle $ABEF$ is convex. In the sequel we will call the configuration of the points $A$ to $J$ thus constructed the \textbf{Harborth configuration}. We will do this, even if the parameter $T$ initially has not been chosen correctly, so that the configuration cannot be completed to the whole Harborth graph.

 Clearly, proceeding as described, with using an arbitrary nonnegative value for the height $T$ of the initial triangle, the line through the final crucial points $H$ and $J$ will only by (a very small) chance meet the line through $A$ and $C$ at an angle $\varphi$ of $90^\circ.$ But, in order that  the Harborth configuration  can be completed to the whole Harborth graph, we have to make sure that this last condition is satisfied.

\begin{figure}
\begin{center}
\includegraphics[width=\linewidth]
{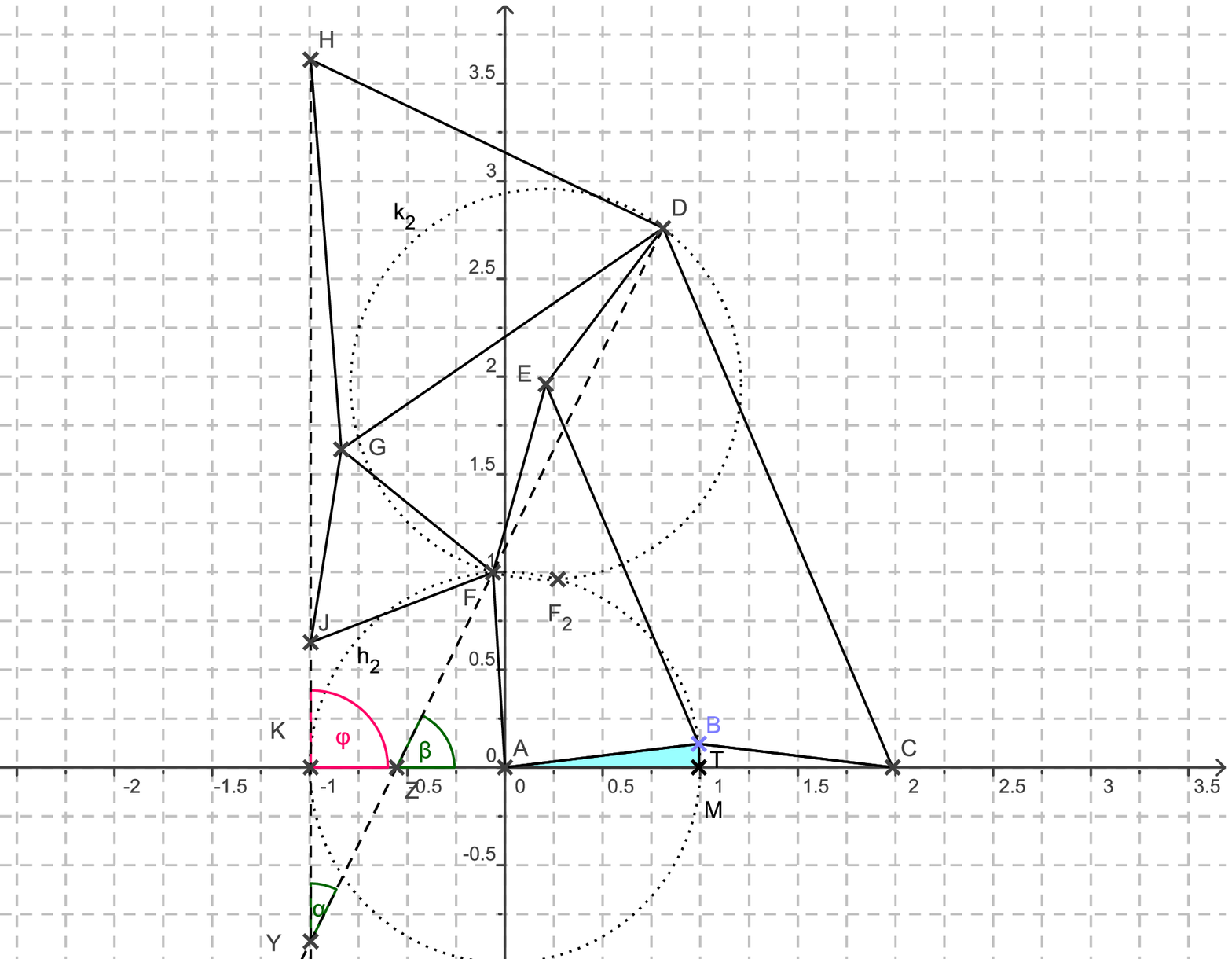}
\caption{A quarter of the Harborth graph, created with GeoGebra \protect\cite{geogebra} - the point $B$ (and thus the height $T$) chosen in order that $\varphi=90^\circ,$ i.e., the configuration can be completed to the whole Harborth graph.}
\end{center}
\label{HarborthPic}
\end{figure}

Using the dynamic geometry software GeoGebra \cite{geogebra}, we are able to manipulate the Harborth configuration in dependence on the parameter $T.$ Furthermore with GeoGebra we are able to read off approximate values (up to five decimal places) for the resulting angles and coordinates. At this early stage of our investigations already, some important observations can be made:

\begin{remark}
\label{1stresults}$\,$\hfill\par
\hangindent
\leftmargini
\textup{(1)}\hskip\labelsep
The point $F,$ and thus the whole Harborth configuration only exists for $T\in[0,b],$ where $b$ is approximately $0.13504.$ \par
\hangindent
\leftmargini
\textup{(2)}\hskip\labelsep
For $T\in [0,b],$ the angle $\varphi$ between $AC$ and $HJ$ lies approximately in the interval $[85.88496^\circ,\,94.59043^\circ].$\par
\hangindent
\leftmargini
\textup{(3)}\hskip\labelsep
For $T\simeq0.12073,$ the angle $\varphi$ approximates $90^\circ.$
\end{remark}

\smallskip
The upper bound $b$ in the above remark is determined by the fact that for $T>b$ the distance  between $A$ and $E$ becomes greater than $2,$ and the circles with centers $A$ and $E$ of radius $1$ do not meet anymore.

As hinted at in the introduction, once we have calculated the exact values for the extremal angles - which we will do in one of the following sections -  the observations collected in Remark \ref{1stresults} allow us to give the first analytic proof\footnote{The basic idea of this proof was first communicated to the author by H.\ L\"owe in 2003. Here it is presented with his kind permission.} of the planarity of the Harborth graph which does not resort to pictures or models only.

\begin{theorem}
\label{TExistenz}
There exists $T\in[0,b]$ such that $\varphi=90^\circ$ exactly, and the Harborth configuration can be completed to the Harborth graph.
\end{theorem}

\begin{proof}{\cite{LoeweHarborth}}
Since the geometric operations we used in our construction (i.e., drawing circles and letting them intersect) depend continuously on their parameters (i.e., centers and radii), and since the composition of continuous functions is again continuous, the angle $\varphi$ depends continuously on the height $T$ as long as the Harborth configuration exists. By the Intermediate Value Theorem the existence of at least one $T$ in the above interval $[0,b]$ is assured such that $\varphi$ is precisely $90^\circ.$
\end{proof}

\begin{corollary}
The Harborth graph indeed is a 4-regular unit-distance graph which is planar. 
\end{corollary}

From now on, our main goal will be to determine a precise description of this particular $T,$ the existence of which we have shown above, without resorting to trial-and-error.

\subsection{Setting up algebraic equations}
\label{settingupCs}  
To describe the points of the Harborth graph more precisely we need to introduce coordinates. For the time being, we choose the point $A$ as the center of our coordinate system, and the ray which extends the base side $AC$ of the triangle $ABC$ as the positive $x$-axis. When considering the Harborth graph as a whole, this is by far not the most obvious choice. In fact, in a later section we will use its center of symmetry as the origin of a more natural coordinate system, and then will have to ''translate'' our intermediate results. As we proceed, we will see that in some places this will prove to be quite cumbersome, which will retroactively justify our initial choice of coordinates. 

Now let $(t,T)$ be the coordinates of the point $B$. Since $A,$ $B,$ and $C$ form an isosceles triangle, the coordinates of the point $C$ are given by $(2t,0).$ 

Denoting the coordinates of any point $P$ of the graph different from $B$ with $(x_P,y_P),$ we have:
\begin{eqnarray}
t^2 + T^2 - 1 & = & 0\label{IniTri}\\
-t\cdot(x_D - 2\cdot t) + T\cdot y_D - \frac{3}{2}& = & 0\label{TrapezD1}\\
(x_D - 2\cdot t)^2 + y_D^2 - 9 & = & 0\label{TrapezD2}\\
t\cdot(x_E - t) - T\cdot(y_E - T) + 1 & = & 0\label{TrapezE1}\\
(x_E - t)^2 + (y_E - T)^2 - 4 & = & 0\label{TrapezE2}\\
x_F^2 + y_F^2 - 1 & = & 0\label{KoordF1}\\
(x_E - x_F)^2 + (y_E - y_F)^2 - 1 & = & 0\label{KoordF2}\\
(x_F - x_G)^2 + (y_F - y_G)^2 - 1 & = & 0\label{KoordG1}\\
(x_D - x_G)^2 + (y_D - y_G)^2 - 4 & = & 0\label{KoordG2}\\
(x_D - x_H)^2 + (y_D - y_H)^2 - 4 & = & 0\label{KoordH1}\\
(x_G - x_H)^2 + (y_G - y_H)^2 - 4 & = & 0\label{KoordH2}\\
(x_F - x_J)^2 + (y_F - y_J)^2 - 1 & = & 0\label{KoordJ1}\\
(x_G - x_J)^2 + (y_G - y_J)^2 - 1 & = & 0\label{KoordJ2}\\
x_H - x_J & = & 0\label{Ortho}
\end{eqnarray}

Let us shortly comment on the meaning of these equations:
(\ref{TrapezD1}), (\ref{TrapezD2}), respectively (\ref{TrapezE1}), (\ref{TrapezE2}), define the vertices $D$ and $E$ of the trapezoid which do not belong to the initial triangle $ABC.$ Equations (\ref{TrapezD1}) and (\ref{TrapezE1}) stem from the fact that the line $BC$ meets $CD$ at an angle of $60^\circ$ and the line $EB$ meets $BC$ at $120^\circ.$ Every further pair of Equations (\ref{KoordF1})-(\ref{KoordJ2}) is chosen in accordance to the geometric constructions described in (\ref{Fconst}) - (\ref{Jconst}), each pair defining one of the points $F$ - $J.$ Finally, Equation (\ref{Ortho}) has to be satisfied in order that the lines given by $AC$ and $HJ$ meet at an angle of $90^\circ.$

Using this set of equations, and approximations for the coordinates, which we read off from Figure 2 with the help of GeoGebra, we are already able to calculate arbitrary precise approximations of all coordinates by using standard numerical algorithms. E.g., the results below show the coordinates with an exactness of $15$ digits; they were calculated with Mathematica, version 4.0.1.0. 
\begin{eqnarray}
B = (t,T)  &\simeq& (0.992685948824186,\, 
0.120725337054926)\nonumber\\
D &\simeq&
(0.809996600722107,\, 
2.760161754567202)\nonumber\\
E &\simeq&
(0.209102417540010,\,
1.960833173433061)\nonumber\\
F &\simeq&
(-0.061398137844065,\,
0.998113354619244)\label{numerics}\\
G &\simeq&
(-0.838419516770942,\,
1.627587561152422)\nonumber\\
H &\simeq&
(-0.995049481192288,\,
3.621444891616507)\nonumber\\
J &\simeq&
(-0.995049481192288,\,
0.639930204451542)\nonumber\\
\nonumber
\end{eqnarray}
\section{Deducing an equation for the $y$-coordinate of the point $B$}

In this section we will deduce the minimal polynomial for the coordinate $y_B=T,$ i.e., the unique primitive integer polynomial $P_T$ of smallest degree such that $P_T(T)=0$ holds \cite[Section 4.1.1]{Cohen}.

Let us shortly describe the main approach which we will repeatedly take: Given two polynomial equations, which are satisfied by certain coordinates, we will calculate the resultant of the corresponding polynomials with respect to one of the appearing variables, thus eliminating this particular variable. Sometimes, if the polynomials are not too complicated, we will use Groebner basis techniques to treat more than two polynomial equations simultaneously. Both procedures lead to new polynomials or polynomial equations, which are consequences of the original ones\footnote{For a more precise description see, e.g., \cite{Cox} or \cite{Loos}.}, but contain fewer variables. To keep expressions from becoming too complicated, and running times from becoming too long we will try to factor the resulting polynomials.  Many times it will prove to be advantageous to allow the factorization to be done over the ring extension $\mathbb Z[\sqrt{3}].$ If a particular polynomial is reducible we will continue our deliberations with that factor, which corresponds to the actual values of the coordinates. To check this, we use numerical approximations analogous to those given in Section \ref{settingupCs}, but which are precise to an error of $\epsilon=10^{-100}$. Most calculations, especially those of resultants, factorizations, and numerical evaluations, were done with Mathematica, version 4.0.1.0.

The succession of eliminations will be determined by the order in which the corresponding points were constructed. E.g., in the initial isosceles triangle $ABC$ due to the choice of the coordinate system all point coordinates are directly expressible in terms of $T.$ Next we will determine polynomials which describe the connection between the coordinates of $D$ and $E,$ respectively, and the parameter $T$. After that the polynomials for the coordinates of $F$ are calculated by using those for the coordinates of $D$ and $E$ and eliminating the variables in between. In principle, continuing this procedure would lead to polynomials in $x_H$ and $T,$ respectively $x_J$ and $T.$ Using the final equation $x_H=x_J,$ one should be able to deduce one polynomial in the variable $T$ alone. Unfortunately, due to the increasing complexity of expressions we were not  able to continue this line of thought to its conclusion, but had to resort to an alternative way. Nevertheless we will try to push as far as possible with this approach, and come up with an alternative, when it proves to be necessary.

In the sequel, we will switch between listing the polynomials and the corresponding polynomial equations at will. When only a polynomial $P$ is given it should be understood that the coordinates appearing in $P$ satisfy the corresponding polynomial equation $P=0.$

\subsection{From $A$ to $F$}

Using the procedure \texttt{GroebnerBasis} with Equations
(\ref{IniTri}) - (\ref{TrapezD2}) as input, like

\begin{verbatim}
GroebnerBasis[{t^2+T^2-1,-t*(xD-2*t)+T*yD-3/2,(xD-2*t)^2+yD^2-9}, 
              {t,xD,yD,T}]
\end{verbatim}

\noindent
one of the polynomials we get is
\begin{equation}
P_{y_D,T}:= 27 - 36 T^2 + 12 T\cdot y_D - 4 y_D^2.
\label{PyDT}
\end{equation}

Analogously, by changing the order of variables,
$$
0=1 - 56 T^2 + 784 T^4 - 8 x_D^2 - 208 T^2\cdot x_D^2 + 16 x_D^4
$$
can be deduced. This last equation is irreducible over $\mathbb{Z},$ but factors over $\mathbb{Z}[\sqrt{3}]$ into
polynomials, which are quadratic in $x_D:$
\begin{equation}
(-1 + 28 T^2 - 12\sqrt{3}T x_D + 4 x_D^2)\cdot
(-1 + 28 T^2 + 12\sqrt{3}T x_D + 4 x_D^2).
\label{xD_T_factored}
\end{equation}

Using numerical results for $T$ and $x_D$ in analogy to (\ref{numerics}), we see that only the first of these polynomials 
\begin{equation}
P_{x_D,T}:= -1 + 28 T^2 - 12\sqrt{3}T x_D + 4 x_D^2
\label{PxDT_redux}
\end{equation}

leads to the correct result. Solving (\ref{PxDT_redux}) and (\ref{PyDT}) for $x_D$ and $y_D,$ respectively, and again discarding those solutions which do not describe the correct coordinates, we get explicit descriptions for the coordinates of $D$ in terms of the parameter $T$: 
\begin{eqnarray}
x_D & = & \frac{1}{2}\left(3\sqrt{3}T+\sqrt{1-T^2}\right),\label{xDT}\\
y_D & = & \frac{3}{2}\left(T+\sqrt{3}\sqrt{1-T^2}\right).\label{yDT}
\end{eqnarray}

Starting with Equations
(\ref{IniTri}),(\ref{TrapezE1}) and (\ref{TrapezE2})
and proceeding in the same manner as above, we are led to those polynomials which describe the dependence of the coordinates of the point $E$ on the parameter $T$:
\begin{equation}
3 T^2 - x_E^2\label{PxET}
\end{equation}
and
\begin{equation}
P_{y_E,T}:= -3 + 7T^2-4Ty_E+y_E^2.\label{PyET}
\end{equation}
Again, the corresponding equations can be explicitly solved:
\begin{eqnarray}
x_E &=& \sqrt{3} T,\label{xET}\\
y_E & = & 2T + \sqrt{3}\sqrt{1-T^2}.\label{yET}
\end{eqnarray}

Next we continue by calculating the coordinates for the point $F$, once more using Mathematica's  \texttt{GroebnerBasis} function. This time we start with the newly found set of Equations (\ref{PxET}) and (\ref{PyET}), together with the Equations (\ref{KoordF1}), (\ref{KoordF2}) defining $F.$ From this we get the following polynomials which describe the dependence of the coordinates $x_F$ and $y_F$ of $F$ on the parameter $T:$
\begin{equation}
\label{PxFTvl}
\begin{minipage}{0.88\linewidth}
{\small
\noindent
$
- 81 + 10800 T^2 - 422496 T^4 + 4272384 T^6 - 19194112 T^8 + 45801472 T^{10} -\\ 
\phantom{-} 63111168 T^{12}
+ 48234496 T^{14} - 16777216 T^{16} +\\
\left(1296 - 92448 T^2 + 1645056 T^4 - 9573888 T^6 
+ 30072832 T^8 -  57655296 T^{10} +\right.\\ 
\left.\phantom{(}  66060288 T^{12} - 29360128 T^{14}\right) x_F^2 + \left(- 7776 + 228096 T^2
- 1555200 T^4 +\right.\\
\left.\phantom{(} 5271552 T^6  - 12189696 T^8 + 15728640 T^{10}- 9437184 T^{12}\right) x_F^4 + \\
\phantom{} \left(20736  - 152064 T^2 + 331776 T^4 - 196608 T^6 - 1310720 T^8 + 2097152 T^{10}\right) x_F^6 +\\
\phantom{} \left(- 20736 + 110592 T^2 - 442368 T^4 + 786432 T^6 - 1048576 T^8\right)x_F^8,
$
}
\end{minipage}
\end{equation}
\begin{equation}
\label{PyFT}
\begin{minipage}{0.8\linewidth}
\noindent
$P_{y_F,T}:=$
{\small
$
81 - 648 T^2 + 144 T^4 - 2304 T^6 + 4096 T^8 + \\  
\phantom{P_{y_F,T}+} \left(432 T - 864 T^3 + 6528 T^5 - 10240 T^7\right) y_F + \\
\phantom{P_{y_F,T}+} \left(- 216 + 1584 T^2 - 5376 T^4 + 9216 T^6\right) y_F^2 + \\ 
\phantom{P_{y_F,T}+} \left( - 576 T + 1536 T^3 - 4096 T^5\right)y_F^3
+ \left(144 - 384 T^2 + 1024 T^4\right) y_F^4.
$
}
\end{minipage}
\end{equation}

\noindent
Once again, the first of these polynomials factors over ${\mathbb Z}[\sqrt{3}]$ into two polynomials of total degree $8$ and degree $4$ in $x_F$. Using the numerical values for $T$ and $x_F,$ we can deduce that only one of these describes the connection between the variables $T$ and $x_F.$ It is
\begin{equation}
\label{PxFT}
\begin{minipage}{0.89\linewidth}
\noindent
$
P_{x_F,T}:=
$
{\small
$
\phantom{(}-9 + 600 T^2 - 3472 T^4 + 5888 T^6 - 4096 T^8 -\\
\phantom{(-} 8\sqrt{3}\left(9 T + 96 T^3 - 112 T^5 + 256 T^7\right) x_F
+ 8\left(9 - 150 T^2 + 256 T^4 - 640 T^6\right) x_F^2 +\\ 
\phantom{(+} 8\sqrt{3}\left(36 T - 96  T^3 + 256 T^5\right) x_F^3 +16\left(- 9 + 24 T^2 - 64 T^4\right) x_F^4.
$
}
\end{minipage}
\end{equation}

\noindent
As said before, trying to continue like this to calculate polynomials for the remaining points $G, H$ and $J,$ will lead into a deadend, because the resulting equations become too unwieldy to handle, and take too much time to calculate, even with the help of Mathematica. Still, our main goal remains to find one single equation describing the parameter $T$ alone. Consequently we have to take a step back, and use a slightly more indirect approach, which we will describe in the section following the next one.

\subsection{Interlude: Calculating the extremal values for which the Harborth configuration exists}
\label{Interlude}
With Equations (\ref{xET}) and (\ref{yET}) thus available, we are able to calculate the exact maximal value for $T,$ hinted at in Remark \ref{1stresults}. To this end, we first observe that in case of $T$ being maximal the line segments $AF$ and $FE$ together, again form a straight line segment of twice the original length. Therefore, for $T$ maximal, the coordinates of the point $E$ satisfy $x_E^2+y_E^2-4=0.$ This, together with (\ref{xET}) and (\ref{yET}), after some small calculation leads to 
\begin{equation}
64T^4-56T^2+1=0.
\end{equation}
Solving for $T,$ and comparing with the numerical values presented in Remark \ref{1stresults}, gives

\begin{lemma}
The minimal and maximal value for $T$ such that the Harborth configuration exists are $T=0,$ and
\begin{equation}
T=\frac{1}{4}\sqrt{7-3\sqrt{5}},
\end{equation}
respectively.
\end{lemma}

Using basic trigonometry, from this we are able to further deduce exact values for the crucial angle $\varphi$ for extremal $T$:

In case of $T=0,$ the points $A,$ $B,$ and $C$ lie on one line. So do the points $D,$ $E,$ and $F.$ Moreover the intersection point $Z$ of these two lines together with $C$ and $D$ form an equilateral triangle, the sides of which have length $3.$ Angle $\beta$ (see Figure 2) becomes one of the angles of this triangle, and thus is equal to $60^\circ.$  Furthermore the points $D,$ $F,$ and $G$ form an isosceles triangle, with the length of the base side $FG$ being one, and the other length being two. An analysis of the triangle formed by the lines $DH,$ and the prolongations of $DF$ and $HJ,$ which contains the triangle $DFG$ completely, allows us to calculate $\alpha.$ Since $\varphi=\alpha+\beta$ some further calculations show

\begin{corollary}
For $T=0,$ the angle $\varphi$ in the Harborth configuration is the unique solution of
\begin{equation}
\sin(\varphi)=\frac{1}{4}(7+3\sqrt{5})\sqrt{\frac{3}{22+6\sqrt{5}}}
\end{equation}  
in the interval $[0,90^\circ],$ which up to an error of $10^{-15}$ is $85.884964999269942^\circ.$
\end{corollary}

As we have already observed, when $T$ attains its maximal value, the points $A,$ $E,$ and $F$ lie on one line, and form the side of the isosceles triangle $ABE.$ Leaving the details to the reader, again only using basic trigonometry - and Mathematica for the calculation of trigonometric expressions - we are able to show

\begin{corollary}
For $T=\frac{1}{4}\sqrt{7-3\sqrt{5}},$ the angles $\alpha$ and $\beta$ in the Harborth configuration (see Figure 2) are the unique solutions of
\begin{align}
\cos(\beta) & =  \frac{\sqrt{3}}{8}\left(\sqrt{3+\sqrt{5}}-\sqrt{7-\sqrt{5}}\right),\nonumber\\
\intertext{and}
\cos(\alpha) & =  \frac{68+3\sqrt{230+34\sqrt{5}}+9\sqrt{5}\left(8+\sqrt{230+34\sqrt{5}}\right)}{2\left(23+3\sqrt{5}\right)\sqrt{97-3\sqrt{5}+3\sqrt{230+34\sqrt{5}}}}\nonumber
\end{align} 
in the interval $[0,90^\circ].$ 

Since $\varphi=\alpha+\beta$,  this leads to $\varphi \simeq 94.590425288952345^\circ$ up to an error of $10^{-15}.$
\end{corollary}

\subsection{From the points $D$ and $F$ to the points $H$ and $J$}
\label{NeuKoord} 
Now we continue with our task of determining the minimal polynomial for that particular $T,$ for which the Harborth graph exists, i.e., for which $\varphi=90^\circ$ holds. Our trick is, not to calculate the coordinates of the points $H$ and $J$ directly in dependence of the second coordinate $T$ of the point $B,$ but to introduce further variables $X$ and $Y$ which will lead to simpler equations. These new variables themselves will depend on the points $D$ and $F.$

\begin{figure}
\begin{center}
\includegraphics[width=\linewidth]
{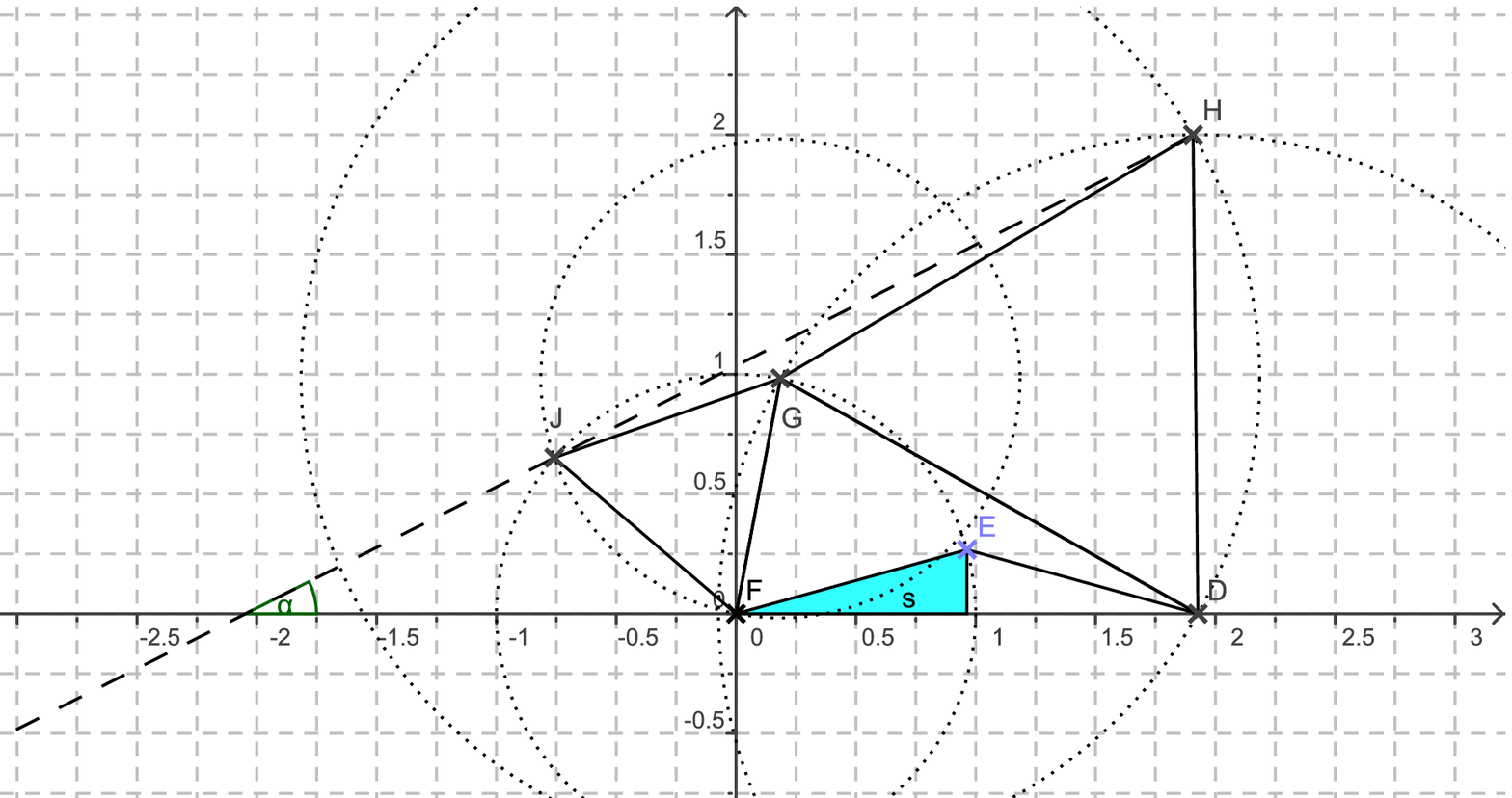}
\caption{Upper part of the Harborth configuration.}
\end{center}
\label{HarborthPic2}
\end{figure}

For this, let us consider $F$ as the origin of a new coordinate system, and the line $FD$ as the new $x$-axis (see Figure 3). Let $(s,S)$ be the coordinates of the point $E$ with respect to this new coordinate system. Clearly, since $DEF$ forms an isosceles triangle, $D$ is described by the coordinates $(2s,0).$ Proceeding as above by using Equations (\ref{KoordG1})-(\ref{KoordJ2}) in this new context, we successively get\footnote{We advise the reader to keep in mind that, although we use the same notation as in Section \ref{settingupCs}, now the coordinates have to be interpreted within the new coordinate frame.}:
\begin{eqnarray}
0&=&3 - 4 s^2 + 4 s x_G,\label{PxGs}\\
0&=&9 - 40 s^2 + 16 s^4 + 16 s^2 y_G^2,\label{PyGs}\\
0&=&9 - 48 s^2 + 48 s^4 + \left(12 s - 48 s^3\right) x_H + 16 s^2 x_H^2,\label{PxHs}\\
0&=&-81 - 144 s^2 - 352 s^4 - 256 s^6 - 256 s^8 +\nonumber\\ 
&& \left(144 s^2 + 896 s^4 + 256 s^6\right) y_H^2 - 256 s^4 y_H^4,\label{PyHs}\\
0&=&9 - 36 s^2 + 16 s^4 + \left(12 s - 16 s^3\right) x_J + 16 s^2 x_J^2,\label{PxJs}\\
0&=&81 - 504 s^2 + 1072 s^4 - 896 s^6 + 256 s^8 +\nonumber\\ 
&& \left(- 144 s^2 + 256 s^4 - 256 s^6\right) y_J^2 + 256 s^4 y_J^4.\label{PyJs}
\end{eqnarray}

All these equations can be easily solved for the respective coordinates. In each case only one of the solutions is in accordance with our geometric construction. Below we present explicit formulas only for the coordinates of the points $H$ and $J$:
\begin{eqnarray}
x_H&=&\frac{-3 + 12 s^2 + \sqrt{3} \sqrt{-9 + 40s^2 - 16 s^4}}{8 s},\\
y_H &=&\frac{3\sqrt{3} + 4\sqrt{3}s^2 + \sqrt{-9 + 40s^2 - 16s^4}}{8 s},\\
x_J&=&\frac{-3 + 4 s^2 - \sqrt{3} \sqrt{-9 + 40 s^2 - 16 s^4}}{8 s},\\
y_J &=&\frac{-3\sqrt{3} + 4\sqrt{3}s^2 + \sqrt{-9 + 40s^2 - 16s^4}}{8s}.
\end{eqnarray}

Therefore the slope of the line $HJ$ with regard to $DF$ as $x$-axis is 
\begin{equation}
\label{firstSteigung}
m_\alpha:=\frac{y_J-y_H}{x_J-x_H}=
\frac{3\sqrt{3}}{4 s^2 + \sqrt{3}\sqrt{-9  + 40 s^2 - 16 s^4}}.
\end{equation}

Now, again we consider the whole Harborth configuration: let new variables $X$ and $Y$ be defined by $X:= x_D-x_F$ and $Y:= y_D-y_F,$ where now $x_D,$ $y_D,$ and $x_F,$ $y_F$ are interpreted as the coordinates of the points $D$ and $F$ with regard to the initial coordinate system. Then the squared length of the line segment $DF$ is given by $X^2+Y^2.$ It follows that 
\begin{equation}
\label{s}
4 s^2=X^2+Y^2,
\end{equation}
 
and Equation (\ref{firstSteigung}) becomes
\begin{equation}
m_\alpha = \frac{3\sqrt{3}}{X^2+Y^2+\sqrt{-9 + 10\left(X^2+Y^2\right) - \left(X^2+Y^2\right)^2}}.
\label{malpha}
\end{equation}

In order to be able to complete the Harborth configuration to the whole Harborth graph, the angle $\varphi$ between the lines $HJ$ and $AC$ must be a right angle. On the other hand, we have $\varphi = \alpha + \beta,$ where $\alpha$ is the angle between $HJ$ and $DF,$ and $\beta$ denotes the angle between $DF$ and $AC,$ as shown in Figure 2. Thus the equality $\alpha = 90^\circ-\beta$ must hold.
Since $0^\circ<\alpha,\beta<90^\circ$, we have $\tan{\alpha}= \tan(90^\circ-\beta)=\frac{1}{\tan(\beta)}.$ Thus the respective slopes satisfy $m_\alpha = 1/m_{\beta}.$ The slope $m_\beta$ of $\beta$ is given within the first coordinate system by $m_\beta= Y/X.$ This, together with (\ref{malpha}), implies
\begin{equation}
\frac{3\sqrt{3}}{X^2+Y^2+\sqrt{-9 + 10\left(X^2+Y^2\right) - \left(X^2+Y^2\right)^2}}=\frac{X}{Y},
\end{equation}

\noindent
which some calculations show to be equivalent to
\begin{equation}
0 = \left(X^2 + Y^2\right)\left(-27 + 30 X^2 - 4 X^4 + 6 \sqrt{3} X Y - 4 X^2 Y^2\right).
\end{equation}
This implies
\begin{equation}
0 = 27 - 30 X^2 + 4 X^4 - 6 \sqrt{3} X Y + 4 X^2 Y^2.
\end{equation}
We set $F(X,Y) := 27 - 30 X^2 + 4 X^4 - 6 \sqrt{3} X Y + 4 X^2 Y^2.$

\noindent
Next we produce polynomials $P_{X,T}$ and $P_{Y,T}$ which describe the connection between $T$ and the new parameters $X,$ and $Y,$ respectively. To do that, this time we use Mathematica's resultant and factorization facilities, as described above. Starting with the polynomial $X-x_D+x_F,$ and the polynomials $P_{x_D,T}$ and $P_{x_F,T}$ given by (\ref{PxDT_redux}) and (\ref{PxFT}), we are thus able to successively eliminate the variables $x_D$ and $x_F$, and get
\begin{equation}
\label{PXT}
\begin{minipage}{0.89\linewidth}
$P_{X,T} :=\\$
{\small
$
\phantom{+\left(\right.}108 T^2 - 684 T^4 + 1344 T^6 - 1344 T^8 +\sqrt{3}\left(- 36 T + 
    300 T^3 - 616 T^5 + 1024 T^7\right) X \\ 
+\left(9 - 213 T^2 + 496 T^4 - 1216 T^6\right) X^2
+\sqrt{3}\left(36 T - 96 T^3 + 256 T^5\right) X^3\\ 
+\left(- 9 + 24 T^2 - 64 T^4\right) X^4\in \mathbb Z[\sqrt{3}][X,T].
$
}
\end{minipage}
\end{equation}

\noindent
In an analogous manner we deduce a polynomial $P_{Y,T}$ in $\mathbb Z[Y,T],$ which for the Harborth configuration describes the connection between these two variables:
\begin{equation}
\label{PYT}
\begin{minipage}{0.88\linewidth}
$P_{Y,T}:=\\$
{\small
$
\phantom{+\left(\right.}
81 - 405 T^2 + 900 T^4 - 1008 T^6 + 448 T^8 
+\left( 54 T - 54 T^3 + 456 T^5 - 512 T^7\right) Y +\\
\phantom{+}\left( - 54 + 207 T^2 - 624 T^4 + 576 T^6\right) Y^2 
+\left(- 18 T + 48 T^3 - 128 T^5\right) Y^3 +\\ 
\phantom{+}\left( 9 - 24 T^2 + 64 T^4\right) Y^4.
$
}
\end{minipage}
\end{equation}

Finally we repeat this procedure with the polynomials $F(X,Y),$ $P_{X,T},$ and $P_{Y,T}$ to eliminate variables $X$ and $T.$ I.e., first we let Mathematica calculate the resultant of $F$ and $P_{X,T}$ with regard to $X.$ We will not present the result here; let it be enough to state that the result is a polynomial of degree $8$ in $Y$, degree $32$ in $T$ and total degree $32,$ which, but for a constant factor, cannot be factored further by Mathematica, even when considered over ${\mathbb Z}[\sqrt{3}].$ With the help of Mathematica, we are able to calculate the resultant of this polynomial and $P_{Y,T}$ with respect to $Y.$ This leaves us with a polynomial in the single variable $T$ of order $156.$ Strangely enough, this final polynomial is reducible over ${\mathbb Z}[\sqrt{3}],$ its factors (up to an integer constant) being $(2T+\sqrt{3}),$ $(2T-\sqrt{3}),$ $(64T^4-24T^2+9)^6,$ and three other integer polynomials of degree $22,$ $28$ and $80,$ respectively. Here, we need only list the polynomial of degree $22,$ since this is the one which has the $y$-coordinate $T$ of the point $B$ of the Harborth graph as one of its real roots\footnote{In fact it is the positive root of smallest modulus.}:
 
\begin{theorem}
The minimal polynomial for the $y$-coordinate $T$ of the vertex $B$ of the Harborth graph is 
\begin{equation}\nonumber
\label{PT}
\begin{minipage}{.97\linewidth}
\noindent
$P_T:=$
{\small
$
-492075 + 52356780 T^2 - 1441635408 T^4 + 12222052416 T^6 
- 60567699456 T^8 +\\
\phantom{P_T=\vert +}  189747007488 T^{10} 
- 417660420096 T^{12} 
+ 607025037312 T^{14} - 655053815808 T^{16} +\\ 
\phantom{P_T=\vert +}  446118756352 T^{18} - 422064422912 T^{20} + 437348466688 T^{22}.
$
}
\end{minipage}
\end{equation}
\end{theorem}

Since this polynomial is primitive, and irreducible over $\mathbb Z,$ it is the minimal polynomial of $y_B=T.$ Thus we have achieved our desired first main result. The rest of this paper is concerned with the determination of the minimal polynomials of the other coordinates, and some of their properties.

\section{Minimal polynomials for the $y$-coordinates}
\subsection{The points $D,E,$ and $F$}
The method to calculate the minimal polynomial for each of the $y$-co\-ordi\-nates of the points $D,E,F$ is very similiar to what we have done in the last subsection. Let $P$ denote one of the points $D,E,$ or $F.$ Using Formulas (\ref{PyDT}), (\ref{PyET}), and (\ref{PyFT}), respectively, which describe the connection between the $y$-coordinate $y_P$ and the parameter $T$ by way of irreducible polynomials $P_{y_P,T}\in{\mathbb Z}[y_P,T]$ of absolute degree $2,2$ and $4,$ we let Mathematica calculate the resultant of $P_{y_P,T}$ and $P_T$ from above with respect to $T$. The result is a polynomial in ${\mathbb Z}[y_P].$ It is of degree $44$ in case of $y_D$ and $y_E$ and of degree $88$ for $y_F$. The first two of these polynomials each split over $\mathbb Z$ into two irreducible factors of degree $22.$ The polynomial for $y_F$ splits into one component of degree $44$ and -- again -- into two polynomials of degree $22.$ In each case one of these factors corresponds to the actual $y$-coordinate. 

\begin{theorem}
The minimal polynomials for $y_D,y_E$ and $y_F$ are

\bigskip
\noindent
$
P_{y_D}:=
$
{\small
$
 -2470693585135788 + 1679453964496051893 y_D^2 - 
    2462573171102886288 y_D^4 +\\ 
\phantom{P_{y_D}=\vert+} 1847147913929328048 y_D^6 - 
    888334179987132288 y_D^8 + 302241307009227264 y_D^{10} -\\ 
\phantom{P_{y_D}=\vert-}  74768143621533696 y_D^{12} + 13516084620361728 y_D^{14} - 
    1721332250836992 y_D^{16} +\\ 
\phantom{P_{y_D}=\vert+} 139442448236544 y_D^{18} - 
    6126808596480 y_D^{20} + 109337116672 y_D^{22},
$
}

\bigskip
\noindent
$
P_{y_E}:=
$
{\small
$
 -387038865725307 + 255845547796716 y_E^2 - 1080696123714384 y_E^4 +\\
\phantom{P_{y_E}=\vert +} 985178573370432 y_E^6 + 290816529555456 y_E^8 + 
    1229422640467968 y_E^{10} -\\
\phantom{P_{y_E}=\vert -} 399291497201664 y_E^{12} - 
    226953868935168 y_E^{14} - 145914316455936 y_E^{16} +\\
\phantom{P_{y_E}=\vert +} 84049703993344 y_E^{18} - 9462031056896 y_E^{20} + 437348466688 y_E^{22},
$
}

\bigskip
\noindent
$
P_{y_F}:=
$
{\small
$
-6156736033068 + 4132620043369020 y_F^2 - 28069535202466347 y_F^4+\\ 
\phantom{P_{y_F}=\vert+} 54174190167055116 y_F^6 - 44321252355544320 y_F^8 + 
    16893977313239424 y_F^{10}-\\ 
\phantom{P_{y_F}=\vert-} 3430375146685440 y_F^{12} + 
    781964817629184 y_F^{14} - 165954075623424 y_F^{16}+\\ 
\phantom{P_{y_F}=\vert+} 16400930701312 y_F^{18} - 579898179584 y_F^{20} + 27334279168 y_F^{22}.
$
}
\end{theorem}

\subsection{The points $G,H,$ and $J$}
\label{yKoordGHJ}
When determining the minimal polynomials of the $y$-coordinates of the points $G,$ $H,$ and $J,$ we have to keep in mind that in Section \ref{NeuKoord} we used a different coordinate system for determining their coordinates than for those of the points $A$ to $F.$ Thus first we will have to transform the former into $y$-coordinates within our original system, which was centered in $A$. Since in this section it is paramount not to confuse these systems, in the sequel we will denote coordinates with regard to the coordinate system centered in $F$ with capital letters, and those with regard to the original one centered in $A$ with small letters. Thus, e.g., Equation (\ref{PyHs}) now becomes
\begin{equation}
\label{PyHsnew}
\nonumber
0=-81 - 144 s^2 - 352 s^4 - 256 s^6 - 256 s^8
+ \left(144 s^2 + 896 s^4 + 256 s^6\right) Y_H^2 - 256 s^4 Y_H^4.
\end{equation}

We will have to treat Equations (\ref{PxGs})-(\ref{PyJs}) accordingly. With this notation, our goal has become to calculate the characteristic polynomials for $y_G,$ $y_H$ and $y_J.$  Elementary analytic geometry tells us that the connection between ''old'' coordinates $(x_P,y_P)$ and ''new'' coordinates $(X_P,Y_P)$ can be described by
\begin{equation}
\label{KoordTrafo}
\begin{pmatrix}
x_P\\ 
y_P
\end{pmatrix}
=
\begin{pmatrix}
x_F\\ 
y_F
\end{pmatrix}
+\frac{X_P}{2s}
\begin{pmatrix}
X\\ 
Y
\end{pmatrix}
+\frac{Y_P}{2s}
\begin{pmatrix}
-Y\\ 
X
\end{pmatrix}
\end{equation}
for $P\in\{G,H,J\},$ where $X=x_D-x_F$ and $Y=y_D-y_F,$ as defined in Section \ref{NeuKoord}. Thus,
\begin{equation}
2s\left(y_P-y_F\right) = X_P\cdot Y + Y_P\cdot X.
\end{equation}
Setting $z_P:=y_P-y_F$ for the moment, successively using Equations (\ref{PxGs})-(\ref{PyJs}), which described the connection between new coordinates of the points $G,$ $H,$ $J,$ and the parameter $s,$ and the equality $4s^2=X^2+Y^2,$ after some calculations in the usual manner, we are able to deduce irreducible polynomials in $\mathbb Z[\sqrt{3}][z_P,X,Y]$ for $P\in\{G,H,J\}.$ Using Equations (\ref{PXT}) and (\ref{PYT}) we are further able to eliminate $X$ and $Y,$ and deduce polynomials in $\mathbb Z[z_P,T]$ of total degree $188$ for each of the points. Each of these splits again, leaving us with irreducible polynomials which are of degree $20$ in $T$ and total degree $20$ for the points $G,J$ and degree $24$ in $T$ and total degree $24$ for the point $H.$ Resubstituting $y_P-y_F$ for $z_P$ and using (\ref{PyFT}) to eliminate $y_F,$ we are left with integer polynomials in variables $y_P$ and $T$ of degree $112$ for the points $G,J$ and one of degree $128$ for $H$, which this time split off irreducible polynomials of total degree $20$ (for $G$ and $H$) and $24$ (for $J$), respectively. Calculating the resultant of these polynomials and $P_T,$ thereby eliminating $T,$ once more in each case we get polynomials of degree $176$ for $y_G,$ $y_H,$ and $y_J,$ respectively. Each contains among others an irreducible factor of degree $22$ -- the minimal polynomial, which we were looking for. Therefore we have:  

\begin{theorem}
The minimal polynomials for the $y$-coordinates of the vertices $G,H, $ and $J$ of the Harborth graph are

\bigskip
\noindent
$
P_{y_G}:=
$
{\small
$
-912811377667500 + 16117998953248125 y_G^2 - 
    36709013218422600 y_G^4 +\\ 
\phantom{P_{y_G}=\vert+} 37940201286814800 y_G^6 - 
    23463887481854208 y_G^8 + 10021184125203456 y_G^{10}-\\ 
\phantom{P_{y_G}=\vert-} 3290335763447808 y_G^{12} + 888521341648896 y_G^{14} - 
    192809455583232 y_G^{16}+\\
\phantom{P_{y_G}=\vert+} 29839017902080 y_G^{18} - 
    2742026240000 y_G^{20} + 109337116672 y_G^{22},
$
}

\bigskip
\noindent
$
P_{y_H}:=
$
{\small
$
-12148787578527675 - 123412000423046805 y_H^2 - 
      441020584930952232 y_H^4+\\ 
\phantom{P_{y_H}=\vert+} 273168911377174014 y_H^6 - 
      27343071784237320 y_H^8 - 3667116898760364 y_H^{10}+\\ 
\phantom{P_{y_H}=\vert+} 823044986987616 y_H^{12} - 32095868573376 y_H^{14} - 
      4779985142784 y_H^{16}+\\ 
\phantom{P_{y_H}=\vert+}  615643279360 y_H^{18} - 27098808320 y_H^{20} + 427098112 y_H^{22},
$
}

\bigskip\noindent
and

\bigskip
\noindent
$
P_{y_J}:=
$
{\small
$
-9964518750000 + 570277711828125 y_J^2 - 1780552966387500 y_J^4+\\ 
\phantom{P_{y_J}=\vert+} 849106838377800 y_J^6 + 644904447905880 y_J^8 - 
    102048280254828 y_J^{10}- \\ 
\phantom{P_{y_J}=\vert-} 56106534718368 y_J^{12}+
    9027433758528 y_J^{14} + 605520976896 y_J^{16}- 103349145600 y_J^{18} -\\ 
\phantom{P_{y_J}=\vert+}  
    2815229952 y_J^{20}+ 427098112 y_J^{22}.
$
}
\end{theorem}

\section{Minimal polynomials for the $x$-coordinates}
As we initially announced, we want to give minimal polynomials for all the coordinates of the most important vertices of the Harborth graph, where the center of origin is supposed to be the center of symmetry of the whole graph, i.e., the point $K$ in Figure 2,  and the axes are the axes of symmetry of the Harborth graph. This means that while the $y$-coordinates remain the same, when we shift our origin from $A$ to $K,$ the $x$-coordinates (with respect to the coordinate system centered in $A$) have to be shifted by $-x_J,$ i.e., we have to set
\begin{equation}
\label{Trafo}
x_P^{\hbox{\scriptsize new}}=x_P - x_J
\end{equation}  
for all points $P,$ where the coordinates on the right hand side denote coordinates with respect to the origin $A.$

As a further difficulty, again we have to pay attention that up above the $x$-coordinates for the points $G,$ $H$ and $J$ were given with respect to yet another, third coordinate system, which had $F$ as its origin and was rotated when considered within the other two coordinate frames.

In the sequel, to avoid misunderstandings, we will switch notation, and denote all $x$-coordinates with respect to the system centered in $A$ by $x_P^{\hbox{\scriptsize old}},$ those with respect to the system centered in $K$ will now become $x_P.$

\subsection{Coordinate transformations for the points $A$ to $F$}
Proceeding step by step as in Section \ref{yKoordGHJ}, but starting from the second equation resulting from (\ref{KoordTrafo}), i.e.,
\begin{equation}
2s (x_J^{\hbox{\scriptsize old}}-x_F^{\hbox{\scriptsize old}}) = X_J\cdot X - Y_J\cdot Y,
\end{equation}

we are able to deduce an irreducible polynomial $P_{{x_J^{\hbox{\scriptsize old}}},T}$ in $\mathbb Z[\sqrt{3}][x_J^{\hbox{\scriptsize old}},T]$ of total degree $24,$ and finally the characteristic polynomial of $x_J^{\hbox{\scriptsize old}}$ of degree $22.$ Since it is an even polynomial\footnote{We call a polynomial \textsl{even} if it only contains monomials of even degree.}, it describes the coordinate $x_A=-x_J^{\hbox{\scriptsize old}}$ as well. Thus,

\begin{theorem}
The minimal polynomial for the $x$-coordinate of the vertex $A$ of the Harborth graph, where the coordinate system is the one given in Figure 1, is

\medskip
\noindent
$
P_{x_A}:=
$
{\small
$
-830376562500 + 1358127000000 x_A^2 - 34144387143750 x_A^4 + 
    96857243056800 x_A^6 -\\ 
\phantom{P_{x_A}=\vert-} 68697978132015 x_A^8 - 189712941147 x_A^{10} + 
    6188723588664 x_A^{12} - 704220643376 x_A^{14} -\\ 
\phantom{P_{x_A}=\vert-} 52577813248 x_A^{16} + 
    27196394496 x_A^{18} - 2918612992 x_A^{20} + 106774528 x_A^{22}.
$
}
\end{theorem}

With the help of the polynomial $P_{x_J^{\hbox{\scriptsize old}},T}$ from above, we can produce a polynomial $P_{\hbox{\scriptsize trafo}}\in\mathbb Z[x_P,x_P^{\hbox{\scriptsize old}},T]$ that describes the coordinate transformation (\ref{Trafo}), by calculating the resultant of $P_{x_J^{\hbox{\scriptsize old}},T}$ and the polynomial $x_P-x_P^{\hbox{\scriptsize old}}+x_J^{\hbox{\scriptsize old}}$ with respect to $x_J^{\hbox{\scriptsize old}}.$ This is of total degree $24,$ of degree $24$ in $T,$ and of degree $8$ in both variables $x_P$ and $x_P^{\hbox{\scriptsize old}}.$ 

Applying the same method as above, i.e., first calculating the resultant of this ''transformation polynomial'' and  the respective polynomials $P_{x_P,T},$ which we now have to interpret as polynomials in the variables $x_P^{\hbox{\scriptsize old}}$ and $T,$ followed by a factorization over $\mathbb Z[\sqrt{3}],$ and finally repeating this process with the resulting polynomial and $P_T$, we are able to deduce the characteristic polynomials of the $x$-coordinates $x_B$ to $x_F.$ Since this procedure should be standard to the reader by now, we will not go into details anymore, but will only present the results. Once again we stress the fact that these polynomials are for $x$-coordinates with respect to $K$ as origin: 

\begin{theorem}
The minimal polynomials for the $x$-coordinates of the vertices of the Harborth graph, where the coordinate system is the one given in Figure 1, are

\medskip
\noindent
$
P_{x_B}:=
$
{\small
$
-17372788157292129 + 85946816541669534 x_B^2 - 
    172967171143553289 x_B^4+\\
\phantom{P_{x_B}=\vert +}  125428630440736260 x_B^6 - 
    35361034276033728 x_B^8 + 4402034757921792 x_B^{10}-\\
\phantom{P_{x_B}=\vert -} 436015591392256 x_B^{12} + 77220067192832 x_B^{14} - 
    11054716223488 x_B^{16}+\\
\phantom{P_{x_B}=\vert +}  874491412480 x_B^{18} - 34734080000 x_B^{20} + 
    557842432 x_B^{22},
$
}

\medskip
\noindent
$
P_{x_C}:=
$
{\small
$
-55268097000787592100 + 83653148035178006805 x_C^2 -\\
\phantom{P_{x_C}=\vert -}  
    49933201015710366166 x_C^4+ 15170804748275250138 x_C^6 -\\ 
\phantom{P_{x_C}=\vert -} 
    2623723693990622868 x_C^8 + 292733387369474292 x_C^{10}-  24051159678783648 x_C^{12} +\\
\phantom{P_{x_C}=\vert +}  1563610131071808 x_C^{14} - 
    77064294460416 x_C^{16}+ 2572257472512 x_C^{18} -\\
\phantom{P_{x_C}=\vert -}  50083921920 x_C^{20} + 427098112 x_C^{22},
$
}

\medskip
\noindent
$
P_{x_D}:=
$
{\small
$
-15937557042969 + 69169635141939 x_D^2 - 133600085051911 x_D^4 +\\ 
\phantom{P_{x_D}=\vert +} 
    150590940104181 x_D^6 - 109441808559384 x_D^8 + 
    53597367271968 x_D^{10} - \\
\phantom{P_{x_C}=\vert -}  17996039805696 x_D^{12} + 
    4144963934208 x_D^{14} - 647005151232 x_D^{16} +\\
\phantom{P_{x_C}=\vert +}  66726690816 x_D^{18} - 
    4293132288 x_D^{20} + 139460608 x_D^{22},
$
}

\medskip
\noindent
$
P_{x_E}:=
$
{\small
$
-30534686672400 + 184473995962680 x_E^2 - 493600710483009 x_E^4 +\\
\phantom{P_{x_E}=\vert +}  
    800738068318020 x_E^6 - 883225203916608 x_E^8 + 
    687262746783744 x_E^{10} - \\
\phantom{P_{x_E}=\vert -} 378024688788480 x_E^{12} + 145061641105408 x_E^{14} - 37695035736064 x_E^{16} + \\
\phantom{P_{x_C}=\vert +} 
    6218402758656 x_E^{18} - 582162055168 x_E^{20} + 27334279168 x_E^{22},
$
}

\medskip
\noindent
$
P_{x_F}:=
$
{\small
$
-622521 + 20028276 x_F^2 - 150285424 x_F^4 - 
    349270464 x_F^6 + 7694997504 x_F^8 -\\
\phantom{P_{x_F}=\vert -} 5213620224 x_F^{10} - 
    109200064512 x_F^{12} + 709185896448 x_F^{14} - 
    1112735219712 x_F^{16}+\\
\phantom{P_{x_F}=\vert +}  387346071552 x_F^{18} + 
    124822487040 x_F^{20} + 8925478912 x_F^{22}.
$
}
\end{theorem}

\subsection{The ''Coup de gr\^ace'' - the minimal polynomial for $x_G$}

Due to the complexity of expressions we were not able to use the above procedure to calculate the characteristic polynomial for $x_G$ - even with the help of Mathematica. Thus we have to resort to one final trick - yet another coordinate system. For this, first we mirror the Harborth configuration. After that we choose the point $J$ as the new origin, and the ray $JH$ as the positive part of the new $x$-axis (see Figure 4).
Consequently, the new $y$-coordinates are the $x$-coordinates from above.

\begin{figure}
\begin{center}
\includegraphics[width=\linewidth]
{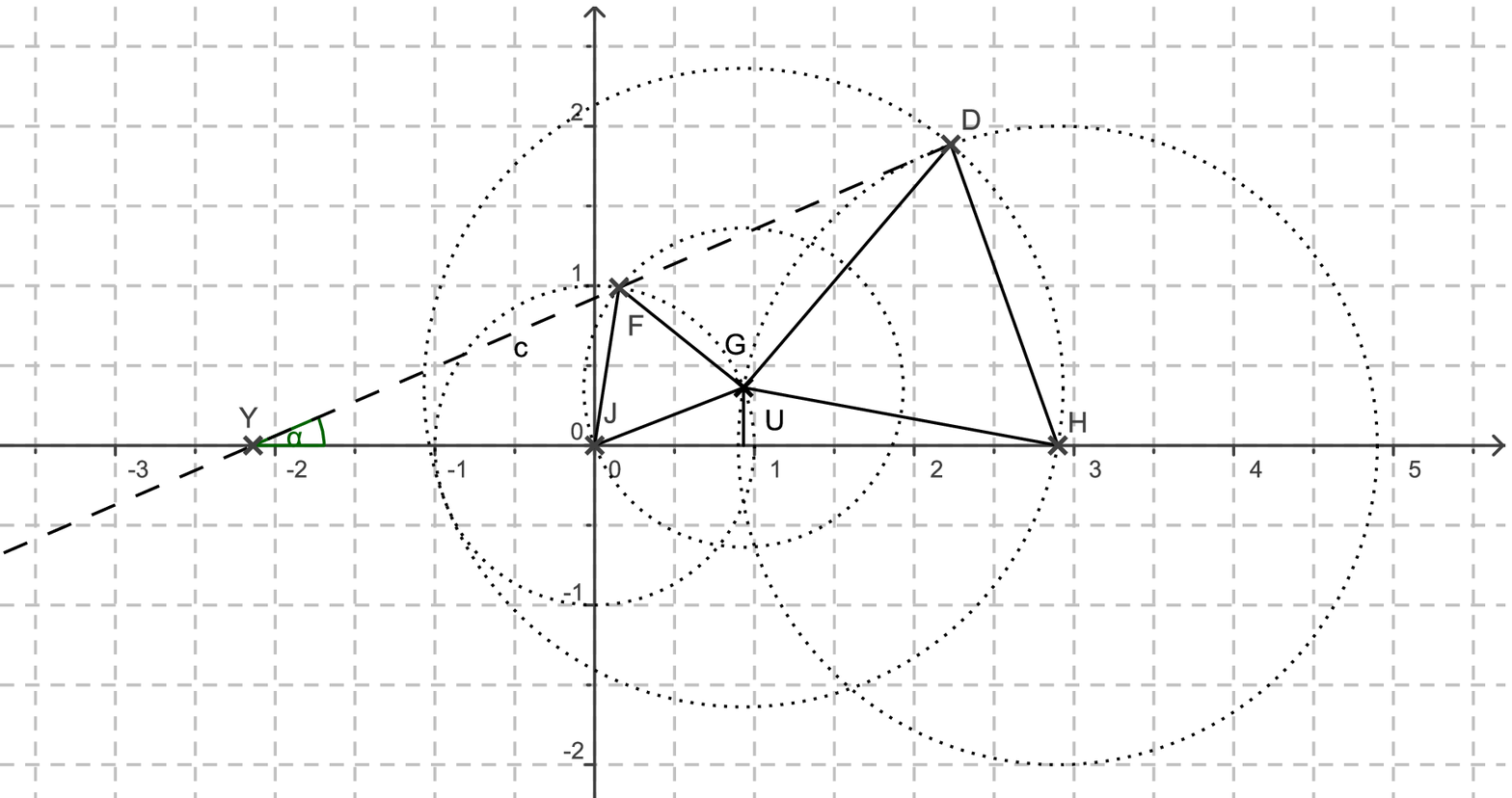}
\caption{Part of the mirrored Harborth graph for the determination of $x_G.$ Note that now $H$ and $J$ lie on the $x$-axis}
\end{center}
\label{Harborth_fuerG}
\end{figure}

Let $(u,U)$ be the coordinates of the point $G$ with respect to this coordinate system. Using Equations (\ref{KoordF2}) and (\ref{KoordG1}), adapted to the new coordinate frame, we get
\begin{equation}
\label{PyFU}
-3+4U^2-4U y_F+4y_F^2=0,
\end{equation}
where $y_F$ now denotes the $y$-coordinate of $F$ with regard to this new system. Since $y_F$ is equal to $x_F$ of old, we thus have calculated a polynomial $P_{x_F,U}.$ Calculating the resultant of this polynomial and $P_{x_F}$ with respect to $x_F$ leads to a polynomial in $\mathbb Z[U]$ of degree $44.$ Factoring this, we achieve an irreducible polynomial of degree $22,$ which is the minimal polynomial of $U$ and thus of $x_G.$ 

\begin{theorem}
The minimal polynomial of the $x$-coordinate of the vertex $G$ of the Harborth graph is

\medskip
\noindent
\begin{center}
\begin{minipage}{.9\linewidth}
$
P_{x_G}:=
$
{\small
$
-106929 + 9380331 x_G^2 - 257190919 x_G^4 + 2410771629 x_G^6 - 
    11872837680 x_G^8+\\
\phantom{P_{x_G}=\vert +}  35430882432 x_G^{10} - 66974055936 x_G^{12} + 
    79549160448 x_G^{14} - 56180293632 x_G^{16}+\\ 
\phantom{P_{x_G}=\vert +} 20514865152 x_G^{18} - 
    2573205504 x_G^{20} + 139460608 x_G^{22}.
$
}
\end{minipage}
\end{center}
\end{theorem}

Thus we have deduced minimal polynomials for all coordinates of the initial Harborth configuration, and consequently, because of the twofold symmetry and our particular, final choice of coordinates, for nearly all the vertices of the Harborth graph. This finishes our initial task. 

\section{Coda}
With all the minimal polynomials at our disposal, we can finish this paper with some observations about their algebraic structure, and consequences for the Harborth graph.

First of all, starting with the minimal polynomials, using a computer algebra system we can once again show the existence of the Harborth graph -- this time by algebraic means only:

\begin{theorem}
For $P\in\{A,\dots,J\}$ let the coordinates $(\pm x_P,\pm y_P)$ of points in the Euclidean plane be particularly chosen roots of the irreducible polynomials $P_{x_P}, P_{y_P}$ which were detailed in the previous sections.
Then these coordinates satisfy the defining equations of coordinates of vertices of the Harborth graph.
\end{theorem}

\begin{proof}
Since the actual calculations have to be done by the computer algebra system and cannot be presented here, the proof will just consist of a series of comments: First, with $(\pm x_P,\pm y_P)$ we denote all four possible combinations of signs of the coordinates, taking into account the twofold symmetry of the Harborth graph, as well as the fact that the polynomials $P_{x_P}$ and $P_{y_P}$ are even. 

Second: Clearly Equations (\ref{IniTri}) - (\ref{Ortho}) have to be restated in accordance with our finally chosen coordinate system, which had the point $K$ as its center. Thus, e.g., instead of (\ref{IniTri}) one has to show that the coordinates of the vertices $A$ and $B$ in the first quadrant satisfy $(x_B-x_A)^2+y_B^2=1,$ i.e., the points are at distance $1$ from each other. Numeric approximations which determine our choice of roots have to be adapted as well.

Moreover calculations with the chosen coordinates are done in the sense of \cite{Loos}, whereat the computer algebra system has to take the full brunt of the work. That is: since there are no explicit formulas for the particular roots of the polynomials $P_{x_P}, P_{y_P},$ which we use as coordinates, we have to understand these roots as completely defined by their representing polynomials, and an isolating interval for each. Instead of an isolating interval one can equivalently use a numeric approximation of sufficient precision. In fact, this seems to be what the computer algebra system Mathematica does: it renders possible computations with algebraic numbers by way of  \texttt{root} objects \cite{MathematicaAlgNumb}. There, choosing a particular root actually means choosing one of the \texttt{root} objects produced by the \texttt{Solve} routine. This we did in accordance with the previously attained numerical results. The ``proof'' of the above theorem has been facilitated then by use of Mathematica's inbuilt \texttt{RootReduce} routine, which had to be applied to the defining equations of the Harborth configuration. In order to achieve that calculations are finished within an acceptable time, the equations have to be expanded and written as sums of products. More concretely, in the example above, after having assigned particular \texttt{root} objects to the variables \texttt{xA}, \texttt{xB}, and \texttt{yB} with the help of the respective polynomials, one successively has to calculate \verb/RootReduce[xA^2]/, \verb/RootReduce[-2*xA*xB]/, etc.\ until one is able to put everything together, thus to be able to calculate \verb/RootReduce[xA^2-2*xA*xB+xB^2+yB^2]/, which indeed results in $1$.
\end{proof}

Even though we have concentrated on one particular set of coordinates, our calculations have shown that, e.g., the $y$-coordinate $y_B$ of the vertex $B$ in any case must satisfy a polynomial equation. Thus there are only finitely many possibilities for coordinates of $B.$ Since by the geometric construction, which we used, all other vertices can be shown to depend uniquely on this one initial vertex and its embedding in the Euclidean plane, different embeddings of the Harborth graph in the plane -- if they were to exist at all -- cannot be transformed into each other in a continuous way. In other words: 
 
\begin{theorem}
The Harborth graph is rigid.
\end{theorem}  

Finally, closer scrutiny of the minimal polynomials leads to

\begin{lemma}
Let $z_P$ be one of the coordinates of a vertex $P$ of the Harborth graph different from zero. Then its minimal polynomial $P_{z_P}\in \mathbb{Z}[X]$ is an even polynomial of degree $22$ and signature $(6,8).$ That is, it has $6$ real zeros and $8$ distinct pairs of conjugate complex zeros.
Consequently we have,
\begin{equation}
\label{radicals}
P_{z_P} = F_{z_P}\circ G
\end{equation}
where $G = X^2\in \mathbb{Z}[X]$, $\circ$ denotes composition and $F_{z_P}$ is a uniquely determined irreducible integer polynomial of degree $11$ and signature $(3,4)$. 
\end{lemma}

A precise proof for this lemma can be done by using any one of the existing algorithms for root isolation and some basic calculus. 

\begin{theorem}
The coordinates of vertices of the Harborth graph which are different from zero cannot be expressed in terms of radicals. Furthermore the Harborth graph as a whole cannot be constructed by compass and ruler alone.
\end{theorem}

\begin{proof}
A well known corollary in Galois theory\footnote{See the corollary to Theorem 46 in \cite{Artin}.} tells us, that the zeros of a polynomial of odd prime degree, which is irreducible over a real number field, are expressible in terms of radicals, if and only if either the polynomial has only one real zero, or all of its zeros are real. We have already ascertained that the polynomials $F_{z_P}$ appearing in Equation (\ref{radicals}) each have $3$ distinct real roots. Thus the equations $F_{z_P}=0,$ and consequently $P_{z_P}=0$ are not soluble by radicals. This  proves the first assumption. 

We show the second assumption by an indirect proof: suppose that the Harborth graph were constructible by compass and ruler. This would imply that the vertex $A$, thus its $x$-coordinate $x_A$ and furthermore $x_A^2$ would be constructible by compass and ruler. As a consequence of this, the order of the Galois group of the corresponding minimal polynomial, which is $F_{x_A}$, would be a power of $2$, and the Galois group would be soluble\footnote{See Theorem 47 in \cite{Artin}.}. This, as we have just seen above, is absurd, completing the proof of the theorem.   
\end{proof}

\bibliographystyle{amsplain}
\bibliography{HarborthV3}
\end{document}